\documentclass[11pt,fleqn,twoside]{article}
\usepackage{amsfonts,amssymb,latexsym}
\makeatletter
\newcommand{\prava}[1]{\small\it
\begin{flushleft}
Copyright \copyright \ 1999 by  #1
\end{flushleft}}

\newcommand{\name}[1]{\begin{flushleft}
                       \LARGE \bf #1
                       \end{flushleft}\vspace{-3mm}}

\newcommand{\Author}[1]{\begin{flushleft}
                       \it #1 \end{flushleft}}

\newcommand{\Adress}[1]{\begin{flushleft}
                       \it #1 \end{flushleft}}

\newcommand{\Date}[1]{\begin{flushleft}
                      \small  \it #1 \end{flushleft}}

\newcommand{\ehkol}{Author \ name}
\newcommand{\ohkol}{Article \ name}
\renewcommand{\@evenhead}{
\hspace*{-3pt}\raisebox{-15pt}[\headheight][0pt]{\vbox{\hbox to \textwidth 
{\thepage \hfil \ehkol}\vskip4pt \hrule}}}
\renewcommand{\@oddhead}{
\hspace*{-3pt}\raisebox{-15pt}[\headheight][0pt]{\vbox{\hbox to \textwidth 
{\ohkol \hfil \thepage}\vskip4pt\hrule}}}
\renewcommand{\@evenfoot}{}
\renewcommand{\@oddfoot}{}

\newcommand{\be}{\begin{equation}}
\newcommand{\ee}{\end{equation}}
\newcommand{\ba}{\hspace*{-5pt}\begin{array}}
\newcommand{\ea}{\end{array}}
\newcommand{\p}{\partial}
\newcommand{\ds}{\displaystyle}
\makeatother

\begin{document}
\setcounter{page}{13}
\thispagestyle{empty}

\renewcommand{\ehkol}{S.F.~Radwan}
\renewcommand{\ohkol}{On the Fourth-Order Accurate  Compact ADI Scheme}

\begin{flushleft}
\footnotesize \sf
Journal of Nonlinear Mathematical Physics \qquad 1999, V.6, N~1,
\pageref{radwa-fp}--\pageref{radwa-lp}.
\hfill {\sc Article}
\end{flushleft}

\vspace{-5mm}

\renewcommand{\footnoterule}{}
{\renewcommand{\thefootnote}{} \footnote{\prava{S.F.~Radwan}}}

\name{On the Fourth-Order Accurate  Compact ADI Scheme for Solving
the Unsteady Nonlinear Coupled Burgers' Equations}\label{radwa-fp}

\Author{Samir F. RADWAN}

\Adress{Department of Engineering Mathematics and Physics, Faculty of
Engineering,\\
Alexandria University, Alexandria, Egypt}

\Date{Received January 26, 1998; Accepted September 8, 1998}

\begin{abstract}
\noindent
The two-dimensional unsteady coupled Burgers' equations with moderate to
severe gradients, are solved numerically using higher-order accurate f\/inite
dif\/ference schemes; namely the fourth-order accurate compact ADI 
scheme, and the fourth-order accurate Du Fort Frankel scheme. The question 
of numerical stability and convergence are presented. Comparisons are
made between the present schemes in terms of accuracy and computational
ef\/f\/iciency for solving problems with severe internal and boundary
gradients. The present study shows that the fourth-order compact ADI
scheme is stable and ef\/f\/icient. 
\end{abstract}

\renewcommand{\theequation}{\arabic{section}.\arabic{equation}}

\section{Introduction}
The basic f\/low equations describing unsteady transport problems form a
parabolic hyperbolic system of partial dif\/ferential equations. The
interaction of the nonlinear convective terms and the dissipative viscous
(or dispersion) terms in these equations can result in relatively
severe gradients in the solution. Also, the accuracy of the numerical
solution and the computational ef\/f\/iciency are highly dependent on the
numerical methods used to solve this kind of partial dif\/ferential equations.
Standard three-point f\/inite dif\/ference methods of approximating spatial
derivatives may work well for smooth solutions, but they fail when severe 
gradients or discontinuities are present, which are common in the shock wave
problems [1, 2, 3].  Lower-order accurate f\/inite dif\/ference methods, such as
upwinding-type f\/inite dif\/ferences, can be a remedy for the numerical
oscillations and dispersions. However, they have a large amount of
``numerical viscosity'' that smooths the solution in much the same way that
physical viscosity would, but to an extent that is unrealistic by several
orders of magnitude, see Sharif and Busnaina [4].  Standard four-point
f\/inite dif\/ference methods, such as Leonard's method [5], are good in their
higher-order accuracies and in reducing numerical smearing ef\/fects. But,
they are plagued by their generation of spurious oscillations or overshoots
in the neighborhood of discontinuities and lack accuracy, as reported by Liu
{\it et al.} [6], and the present author [7].  TVD f\/inite dif\/ference
schemes~[8, 9], guarantee oscillation-free solutions, but they are limited to
second-order accuracy. Third-order accurate TVD schemes are reported by
Gupta {\it et al.} [10], and by Liu {\it et al.}~[6].  However, in their studies,
nonlinearities were either not present or played a minor role. Therefore,
the perfect numerical methods should possess both higher-order accuracy and
sharp resolution of discontinuities without excessive smearing. Moreover,
higher-order accurate numerical methods are attractive for problems with
long computational time or with required higher accuracy solutions.  As
mentioned by Orszag [11], that at the cost of slight additional
computational complexity, the fourth-order schemes achieve results in the 5\%
accuracy range with approximately half the spatial resolution in each space
direction compared with the second-order schemes (i.e. a factor 8 fewer
grid points in three dimensions).  But, the objection to the standard
higher-order schemes comes from the additional nodes necessary to achieve
the higher-order accuracy.  This precludes the use of implicit methods since
the obtained matrix is not of tridiagonal form, and it is necessary to use
f\/ictitious nodes for the boundary conditions.  Also, they do not allow
easily for non-uniform grids, unless at the expense of the order of
accuracy.  On the other hand, the compact schemes that treat the function
and its necessary derivatives as unknowns at the grid nodes, like the Pade
scheme [40], are fourth-order accurate, and compact in sence that they
reduce to tridiagonal form.  The compact schemes generally consist of f\/inite
dif\/ference schemes which involve two or three grid points.  The three-point
schemes fall into two classes. The f\/irst class consists of methods which are
fourth-order accurate for uniform grids, such as Kreiss scheme [11, 12], the
Mehrstellen method [13], the operator compact implicit scheme [14--16], and
the Hermitian f\/inite dif\/ference method of Peters~[17].  The second class
consists of methods that allow variable grids such as the cubic spline
methods of Rubin and Graves [18--21], and the Hermitian f\/inite dif\/ference
method of Adam [22, 23].  For incompressible viscous f\/low problems, there has
been several work on the construction of compact schemes for the
incompressible Navier-Stokes equations, see for example, [24--33].  The most
noted ones include the work of Gupta [24, 25]. He introduced a compact
fourth-order f\/inite dif\/ference scheme with three nodal points for the
convection dif\/fusion equations. His scheme does not seem to suf\/fer
excessively from spurious oscillatory behavior or numerical viscosity, and
it converges with standard methods such as Gauss Seidel or SOR regardless of
the dif\/fusion [34, 35].  However, Yavneh [36], in his analysis of Gupta's
scheme for the 2-D convection-dif\/fusion equation with Dirichlet boundary
conditions, has concluded that Gupta's compact scheme does not suf\/fer from a
cross-stream artif\/icial viscosity, but it does include a streamwise
artif\/icial viscosity that is inversely proportional to the natural
viscosity.  For unsteady viscous f\/low cases, Weinan and Liu [37] introduced
an modif\/ied version of Gupta's compact scheme based on the vorticity-stream
function formulation. The requirement of compactness of this scheme is
slighly relaxed (for the convection terms) so that the resulting scheme is
simple, easy to implement, and with small phase error.  They were successful
to compute the driven cavity problem at high Reynold numbers upto $10^6$ on
$1024^2$ grid, and they have concluded that the fourth-order compact methods
are comparable in accuracy with the spectral methods for most problems of
practical interests.  For turbulent f\/luid f\/low, where there is a range of
space and time scales, Lele [38] had introduced a series of higher-order
compact schemes that are generalization of the Pade scheme with three nodal
points, and with an improved representation of the shorter length scales.
He applied them to the evolution of  supersonic shear layers.      

The disadvantage of  the above higher-order compact schemes involving
three nodal points is that the boundary conditions are no longer suf\/f\/icient
and they do not allow easily for non-uniform grids, unless at the expense of
the order of accuracy.  Another disadvantage of some compact schemes is the
complexity of the resulting nonlinear f\/inite dif\/ference equations and the
associated dif\/f\/iculty in solving them ef\/f\/iciently.  On the other hand, the
compact scheme with two nodal points, like second-diagonal Pade scheme, is
fourth-order accurate even for non-uniform spatial grids, and no f\/ictitious
points neither extra formula are needed for Dirichlet boundary conditions,
as discussed by White [39] and Keller [40].  Also, Liniger and Willoughby
[41] studied the numerical solutions of stif\/f systems of ordinary
dif\/ferential equations, that are encountered in many areas of applied
mathematics, using compact two-point implicit methods.  They introduced
three main compact schemes with dif\/ferent order of accuracy, and with some
very favorable properties.  In particular, their schemes have $A$-stability
in the sense of Dahlquist [55] and they account for the exponential
character of the rapidly decaying solutions directly, which are referred as
exponential f\/itting methods.  

Inspite of many articles have appeared in the literature concerning the
applications of the higher-order accurate schemes including the compact
schemes to f\/luid dynamics problems, there is no much works done in the area
of application of two-point compact schemes [39--44], like the fourth-order
accurate second-diagonal Pade approximation, to multi-dimensional cases.
This is the main objective of present study, where we study the feasibility
of extending the two-point compact scheme to solve the unsteady
two-dimensional coupled Burgers' equations. They take the following 
form:
\be \label{radwa:1.1}
\frac{\p u}{\p t} =-u\frac{\p u}{\p x} -v\frac{\p u}{\p y}+
\nu\frac{\p^2 u}{\p x^2} +\nu\frac{\p^2 u}{\p y^2},
\ee
\be \label{radwa:1.2}
\frac{\p v}{\p t} =-u\frac{\p v}{\p x} -v\frac{\p v}{\p y}+
\nu\frac{\p^2 v}{\p x^2} +\nu\frac{\p^2 v}{\p y^2},
\ee
\[
x_0\leq x\leq x_N, \qquad y_0 \leq y \leq y_M, \qquad t>0
\]
with the initial conditions
\be \label{radwa:1.3}
u(x,y,0)=u_0(x,y), \qquad v(x,y,0)=v_0(x,y)
\ee
and the Dirichlet boundary conditions
\be \label{radwa:1.4}
\ba{lll}
u(x_0,y,t)=u_1(y), & \qquad & u(x_N,y,t)=u_2(y),
\vspace{2mm}\\
u(x,y_0,t)=u_3(x), & & u(x,y_M,t)=u_4(x),
\ea
\ee
\be \label{radwa:1.5}
\ba{lll}
v(x_0,y,t)=v_1(y), & \qquad  & v(x_N,y,t)=v_2(y),
\vspace{2mm}\\
v(x,y_0,t)=v_3(x), & & v(x,y_M,t)=v_4(x),
\ea
\ee
where $\nu$ is equal to $1/Re > 0$, and $Re$ is the Reynolds
number. Coupled Burgers' equations are used to model many practical
transport problems, such as vorticity transport, hydrodynamic turbulence,
shock wave theory, wave processes in thermoelastic medium, transport and
dispersion of pollutants in rivers, and sediment transport, see references
[45--47].  The coupled Burgers' equations are an appropiate form of the
Navier-Stokes equations. They have the same convective and dif\/fusion form as
the incompressible Navier-Stokes equations. Fletcher [48] made a comparison
of f\/inite element and f\/inite dif\/ference methods with dif\/ferent orders of
accuracy for solving the two-dimensional Burgers' coupled equations, and he
concluded that the f\/ive-point f\/inite
dif\/ference scheme is the most ef\/f\/icient
scheme. Boonkkamp and Verwer [49] have used the extrapolated Odd-Even
Hopscotch scheme for solving the inhomogeneous two-dimensional coupled
Burgers' equations, but their solution exhibited wiggles.  Arminjon and
Beauchamp [50] have concluded that the f\/inite element method is ef\/f\/icient
compared to the other methods, namely the method of lines and
Runge-Kutta-type method in solving the above Burgers' equations.  Jain and
Raja [51] have used splitting-up technique to reduce the problem to a
sequence of tridiagonal systems. Jain and Lohar [52] used spline locally one
dimensional (SOLD) algorithm for solving coupled nonlinear parabolic
equations. EL-Zoheiry and EL-Naggar [53] used the spline alternating
direction implicit (SADI) method for solving the two-dimensional Burgers'
equations.  EL-Naggar [54] has presented a mixed implicit-explicit two
levels algorithm that is based on SADI method for solving the coupled
Burgers' equations. However, his results exhibited oscillations.

In the present study, higher-order accurate two-point compact alternating
direction implicit algorithm (CADI) is introduced to solve the
two-dimensional unsteady coupled Burgers' equations, for problems with
moderate to severe internal and boundary gradients.  The algorithm has the
following features:
\begin{enumerate}
\item[(1)] it results in f\/inite dif\/ference equations that involve only
two-nodal points and therefore is formally fourth-order accurate on all grid
points, even for non-uniform grids;
\item[(2)] it has $A$-stability in the sense of Dahlquist, and accounts for
the exponential character of rapidly varying solutions [41];
\item[(3)] it utilizes Newton's method for linearization with a quardratic
convergence; 
\item[(4)] it requires only the given Dirichlet boundary conditions;
\item[(5)]  the algorithm is simple and easy to implement.
\end{enumerate}
In short, the present method is the natural extension of $A$-stable
fourth-order accurate second-diagonal Pade approximation to solve
multi-dimensional f\/low problems with moderate to severe gradients.
Comparison of the present scheme with the fourth-order Du Fort Frankel
scheme is made in terms of accuracy and computational ef\/f\/iciency, which show
that the fourth-order compact ADI scheme is stable and ef\/f\/icient.

\setcounter{equation}{0}

\section{The numerical schemes}

In this section, the present numerical schemes, namely the fourth-order
accurate two-point compact ADI scheme and the fourth order Du Fort Frankel
scheme, are derived for the two-dimensional unsteady coupled Burgers'
equations (1). 

\subsection{Fourth-order accurate Du Fort Frankel scheme}

Let the interval $[x_0, x_N]$ be discretized into $N$ grid steps of size 
$\Delta x$, where $\Delta x = ( x_i - x_{i-1} )$, $i$ is an index of any
grid-point in $x$ direction. Similarly, the interval $[y_0, y_M]$ is
discretized into $M$ grid steps of size  $\Delta y$,  where  
$\Delta y = (y_j - y_{j-1})$, $j$ is an index of any grid point in
$y$-direction, and $n$ is an index for the temporal grid point. The explicit
form of the Du Fort Frankel scheme for the two-dimensional coupled Burgers'
equations (1), using Kreiss fourth-order accurate approximations [11] for
the spatial derivatives, takes the following form: 
\be \label{radwa:2.1}
\ba{l}
\ds \left[ \frac{\vec \psi_{ij}^{n+1} -\vec \psi_{ij}^{n-1}}{2\Delta t}
\right] =-u_{ij}^n D_x \left[ 1-\frac{\Delta x^2}{6} \delta_x^2 \right]
\vec \psi_{ij}^n + \nu \delta_x^2 \left[ 1-\frac{\Delta x^2}{12}\delta_x^2\right]
\vec \psi_{ij}^n 
\vspace{3mm}\\
\ds \qquad - 
v_{ij}^n D_y \left[ 1-\frac{\Delta y^2}{6} \delta_y^2 \right]
\vec \psi_{ij}^n + \nu \delta_y^2 \left[ 1-\frac{\Delta y^2}{12}
\delta_y^2\right]
\vec \psi_{ij}^n ,
\ea
\ee
where
\be \label{radwa:2.2}
\vec \psi =[u,v]^T,
\ee
\be \label{radwa:2.3}
D_x \vec \psi_{ij} =\frac{1}{2\Delta x} \left(\vec \psi_{i+1j}-
\vec \psi_{i-1j}\right),
\ee
\be \label{radwa:2.4}
\delta_x^2 \vec \psi_{ij} =\frac{1}{\Delta x^2}
\left( \vec \psi_{i+1 j} -2 \vec \psi_{ij} +\vec \psi_{i-1 j}\right).
\ee

Def\/ine  $\ds c_x^n=u_{ij}^n \frac{\Delta t}{\Delta x}$, 
$\ds c_y^n =v_{ij}^n \frac{\Delta t}{\Delta y}$, to be the local Courant
numbers in $x$ and $y$ directions, 
$\ds d_x =\nu\frac{\Delta t}{\Delta x^2}$, and
$\ds d_y =\nu \frac{\Delta t}{\Delta y^2}$.
The above equation represents the fourth-order accurate leap-frog scheme for
equations (1), and in order to obtain the f\/inal form of the fourth-order
accurate explicit Du Fort Frankel scheme for the 2-D unsteady coupled
Burgers' equations~(1), the center node value $(\psi_{ij})$ in the dif\/fusion
terms in equation (2.1) are replaced by their average at time-levels $(n-1)$
and $(n+1)$, giving: 
\be \label{radwa:2.5}
\ba{l}
\ds \vec \psi_{ij}^{n+1} = A \vec \psi_{ij}^{n-1} + B
\vec \psi_{i+2j}^{n} + C \vec \psi_{i+1j}^{n} + D \vec \psi_{i-1j}^{n}
+ E \vec \psi_{i-2j}^{n}
\vspace{1mm}\\
\ds \qquad + F\vec \psi_{ij+2}^{n}+ G \vec \psi_{ij+1}^{n}+
H\vec\psi_{ij-1}^{n} +L \vec \psi_{ij-2}^{n},
\ea
\ee
where
\be \label{radwa:2.6}
\ba{lll}
A=(1-2.5d_x -2.5 d_y)/Q, & \qquad &B=(c_x-d_x)/6Q,
\vspace{1mm}\\
C=(-8c_x^n+16d_x)/6Q, & & D=(8c_x^n+16d_x)/6Q,
\vspace{1mm}\\
E=-(c_x^n+d_x)/6Q, & & F=(c_y^n-d_y)/6Q,
\vspace{1mm}\\
G=(-8c_y^n+16d_y)/6Q, & & H=(8c_y^n+16d_y)/6Q,
\vspace{1mm}\\
L=-(c_y^n+d_y)/6Q, & & Q=(1+2.5d_x+2.5d_y).
\ea
\ee

\subsection{Fourth-order accurate compact ADI scheme}

Liniger {\it et al.} [41] have introduced  the following linear one step
formulas for $\phi(x)$ containing real free parameters ($a$ \& $b$):
\be \label{radwa:2.7}
[\phi_{i+1}-\phi_i]-\frac{\Delta x}{2} [(1+a)\phi_{x_{i+1}}+(1-a)\phi_{x_i}]
+\frac{\Delta x^2}{4} [(b+a) \phi_{xx_{i+1}} -(b-a)\phi_{xx_i}]=e_T,
\ee
\be \label{radwa:2.8}
e_T =\frac{\Delta x^3}{4} \int_0^1 [ 2\xi^2 -2(1-a)\xi +(b-a)]
\frac{\p^3}{\p x^3} \phi(x+\xi \Delta x) d\xi.
\ee
For the case of $a\geq 0$ and $b=1/3$, the resulting formula has a
third-order accuracy.  Moreover, for the case of $a=0$ and $b=1/3$, the
formula has a fourth-order accuracy, which is known as the two-point
second-diagonal Pade approximation: 
\be \label{radwa:2.9}
[\phi_{i+1} -\phi_i] -\frac{\Delta x}{2} [\phi_{x_{i+1}}+\phi_{x_i}]
+\frac{\Delta x^2}{12} [\Phi_{xx_{i+1}} -\phi_{xx_i} ] =e_T,
\ee
\be \label{radwa:2.10}
e_T =\frac{\Delta x^5}{24} \int_0^1 \xi^2(\xi-1)^2 \frac{\p^5}{\p x^5}
\phi(x+\xi \Delta x) d\xi,
\ee
where $\phi_x$, $\phi_{xx}$ are the f\/irst and the second derivatives of the
function $\phi(x)$. Using the above scheme and an ADI-type time marching
procedure for the temporal derivative, the compact alternating direction
implicit algorithm (CADI) for the coupled Burgers' equations (1), are
obtained by f\/irst rewriting these equations as follows:
\be \label{radwa:2.11}
\left[ \ba{c} u \\ v \ea \right]_t =
\left[ \ba{c} (\nu u_x -0.5 u^2)_x
\vspace{1mm}\\
(\nu v_x -u v)_x +v u_x
\ea \right]+
\left[ \ba{c}
(\nu u_y -uv)_y +uv_y
\vspace{1mm}\\
(\nu v_y -0.5 v^2)_y \ea \right],
\ee
then the ADI -- type time marching procedure requires, in one full time
step, the solution of: 

{\bfseries \itshape $x$-sweep}
\renewcommand{\theequation}{\arabic{section}.\arabic{equation}{\rm a}} 
\setcounter{equation}{11}
\be \label{radwa:2.12a}
\left[ \ba{c} u_t \\ v_t \ea \right]^{n+0.5} =
\left[ \ba{c} (\nu u_x -0.5 u^2)_x
\vspace{1mm}\\
(\nu v_x -u v)_x +v u_x
\ea \right]^{n+0.5}+
\left[ \ba{c}
g_1
\vspace{1mm}\\
g_2
\ea \right]^n,
\ee
\renewcommand{\theequation}{\arabic{section}.\arabic{equation}{\rm b}} 
\setcounter{equation}{11}
\be \label{radwa:2.12b}
\left[ \ba{c} u_x\\ v_x \ea \right]^{n+0.5} =
\left[ \ba{c} F \\ G \ea \right]^{n+0.5},
\ee

{\bfseries \itshape $y$-sweep}
\renewcommand{\theequation}{\arabic{section}.\arabic{equation}{\rm a}} 
\setcounter{equation}{12}
\be \label{radwa:2.13a}
\left[ \ba{c} u_t \\ v_t \ea \right]^{n+1}=
\left[ \ba{c}
(\nu u_y -uv)_y +uv_y
\vspace{1mm}\\
(\nu v_y -0.5 v^2)_y \ea \right]^{n+1}+\left[ \ba{c} f_1 \\ f_2 \ea
\right]^{n+0.5} ,
\ee
\renewcommand{\theequation}{\arabic{section}.\arabic{equation}{\rm b}} 
\setcounter{equation}{12}
\be \label{radwa:2.13b}
\left[ \ba{c} u_y\\ v_y \ea \right]^{n+1} =
\left[ \ba{c} H \\ T\ea \right]^{n+1},
\ee
where, $[f_1, f_2]^T$, $[g_1 , g_2]^T$ are the components of the f\/irst and
the second matrices in the right hand side of equation (2.11). The solution
procedure consists of solving, f\/irst, equations (2.12) in the solution
vector $[U, V, F, G]^T$ at time level $n +0.5$, ($x$-sweep), then solving
equations (2.13) in the solution vector $[U, V, H, T]^T$ at time level
$n+1$, ($y$-sweep).  Noting that $\alpha =1/\Delta t$, and in order to apply
the compact scheme to the solution in the $x$-sweep, a vector $\vec Q$ and
its derivatives with respect to $x$, for Burgers equations (2.12), have been
def\/ined as follows: 
\renewcommand{\theequation}{\arabic{section}.\arabic{equation}} 
\setcounter{equation}{13}
\be \label{radwa:2.14}
\vec Q_{ij}^{n+0.5}= \left[ \ba{c}
\nu u_x -0.5 u^2\\
\nu v_x -uv\\
\nu u\\
\nu v \ea \right]^{n+0.5}_{ij}=
\left[\ba{c}
\nu F -0.5 U^2\\
\nu G -UV\\
\nu U\\
\nu V \ea \right]^{n+0.5}_{ij},
\ee
\be \label{radwa:2.15}
\vec Q_{xij}^{n+0.5} = \left [ \ba{c}
u_t\\
v_t-v u_x\\
\nu F\\
\nu G
\ea \right]_{ij}^{n+0.5}- \left[\ba{c} g_1 \\
g_2 \\ 0 \\ 0 \ea \right]^n_{ij}
= \left[ \ba{c} \alpha U \\ \alpha  V -VF \\ \nu F \\ \nu G \ea
\right]^{n+0.5}_{ij}-
\left[ \ba{c} g_1 +\alpha U\\ g_2 +\alpha V \\ 0 \\ 0 \ea \right]_{ij}^{n},
\ee
\be \label{radwa:2.16}
\ba{l}
\vec Q_{xxij}^{n+0.5} = \left [ \ba{c}
F_t\\
G_t-vF_x-Fv_x\\
u_t +uu_x\\
v_t +u v_x
\ea \right]_{ij}^{n+0.5}-
\left[\ba{c} g_{1_x} \\ g_{2_x} \\ g_1 \\ g_2 \ea \right]^n_{ij}
\vspace{3mm}\\
\ds \qquad
= \left[ \ba{c} \alpha F \\ \alpha  G -VF_x-FG \\
\alpha U+UF \\ \alpha V +UG \ea
\right]^{n+0.5}_{ij}-
\left[ \ba{c}
\alpha F +g_{1_x}\\
\alpha G + g_{2_x} \\ \alpha U +g_1 \\ \alpha V +g_2
 \ea \right]_{ij}^{n},
 \ea
\ee
where $g_1$, $g_2$, $g_{1x}$, $g_{2x}$ and $F_x$ in equations (2.15)--(2.16)
are approximated by fourth-order accurate f\/inite dif\/ferences. Having
substituted the vector $\vec Q$ and its derivatives into the above two-point
second-diagonal Pade approximation, equation (2.9),
by replacing~$\phi_i$ by
the vector $\vec Q$, we have four non-linear coupled f\/inite dif\/ference
equations in the solution vector $[U, V, F, G]^T$. Newton's method is used
to linearize the equations, and the numerical solution is obtained by
iteration.  The resulting linearized equations form a block tridiagonal
matrix system of order $N$, as  in the following form:
\be \label{radwa:2.17}
a_i \vec \delta_{i-1} +b_i \vec \delta_i +c_i \vec \delta_{i+1} =\vec r_i,
\qquad i=1,2,\ldots,N,
\ee
where $a_i$, $b_i$, and $c_i$ are block matrices of order four, 
$\vec \delta =[\delta U, \delta V, \delta F, \delta G]^T$ is the
change in the solution vector, and $\vec r$ is the right handside vector,
each of order four.  At each iteration, the LU-factorization algorithm is
used to obtain the solution of the system~(2.17). Similarly, the solution
procedure of the Burgers' equation (2.11) in the $y$-sweep, using equations
(2.13).

\subsection{Numerical stability limits}

The implicit formulation of the two-point compact scheme to the
Burgers' equations is always unconditionally stable. In this case,
the accuracy of the numerical solution depends on the size of
the discretizations,
and higher accuracy can be obtained by f\/iner discretization. Moreover,
the present higher-order scheme allows us to use large discretization in
comparison with the second-order schemes.  Also, it is well known that, for
the convection dif\/fusion equation, the leap-frog scheme is
unconditionally unstable,  while the Du Fort Frankel scheme has a stability
limit ($c\leq 1$ \& $d>0$)~[56].
Therefore, it is necessary to use Von Neumann
stability analysis to def\/ine the stability limit of the
fourth-order accurate
Du Fort Frankel scheme to Burgers' equations.   In the present study, it
will be suf\/f\/icient to examine the stability of the fourth-order Du Fort
Frankel scheme for one of the above coupled Burgers' equations (2.5), say
for $u$-component. Let the numerical solution $U(x,y,t)$
be represented by a f\/inite Fourier series,  and for linear stability,
we can examine the behaviour of a single term of the series, as follows:
\be \label{radwa:2.18}
u(i\Delta x, j \Delta y, n\Delta t)= G(n\Delta t)e^{I[k_x i \Delta x+k_y j
\Delta y]},
\ee
where $G(t)$ is the amplitude function at time-level n of this term whose
wave numbers in the $x$ and $y$ directions are  $k_x$ and $k_y$,
and $I=\sqrt{-1}$.    Def\/ining the $x$ and $y$  phase angles as
$\theta _x=k_x \Delta x$  and   $\theta_y  =  k_y \Delta y$,
then, equation (2.18) becomes:
\be \label{radwa:2.19}
u^n_{ij}=G^n e^{I[i\theta_x+j\theta_y]}.
\ee
Substituting (2.19) into the f\/irst equation of (2.5), we obtain a
quadratic equation for the amplif\/ication factor $\zeta$,
its solution is:
\be \label{radwa:2.20}
\zeta =\frac{G^{n+1}}{G^n} =\frac 12 \left[ \lambda \pm \sqrt{\lambda^2 +4A}
\right]
\ee
and
\be \label{radwa:2.21}
\ba{l}
\ds \lambda =\frac{1}{3Q}\Bigl\{ \bigl[ (16 \cos \theta_x -\cos 2\theta_x)d_x
+(16\cos \theta_y -\cos 2\theta_y) d_y\bigr]
\vspace{3mm}\\
\ds \qquad + I \bigl[ (-8\sin\theta_x+\sin 2\theta_x)c_x +(-8\sin \theta_y +
\sin 2\theta_y)c_y\bigr]\Bigr\},
\ea
\ee
where $A$ and $Q$ are def\/ined by equations (2.6). For the special case of
$d_x= d_y = d$ and $c_x = c_y = c$, the modulus of the amplif\/ication
factor $\zeta$, def\/ined by the following equation:
\be \label{radwa:2.22}
\chi(c,d,\theta_x,\theta_y)=\max \left( \left| \frac 12
\left[ \lambda +\sqrt{\lambda^2 +4A}\right]\right|, \left|\frac 12
\left[ \lambda -\sqrt{\lambda^2 +4A}\right]\right|\right)
\ee
has been computed for dif\/ferent values of $|c|$ and $d$ and plotted,
as shown in Fig.~1.  This shows that the fourth-order accurate
Du Fort Frankel scheme is unstable for the range
$(0.35\leq |c|\leq 1.0)$.   E.g.   $\chi(1,0.5,\pi/2, \pi/2)=1.77$,
$\chi(1,0.5,\pi, \pi)=1.29$,
$\chi(0.5,0.5,\pi/2$, $\pi/2)=1.14$,
and $\chi(0.5,0.5,\pi, \pi)=1.29$.
For small values of  $|c|$ and $d$  ($|c| < 0.35$), the
instability only occurs for phase angles $\theta_x$ and $\theta_y$
close to $\pi$.  Moreover,  for smaller  values  of  $d$
($d< 0.1$), the  scheme  has a neutral  stability ($\chi =1.0$).
 E.g. $\chi(0.25,0.5,\pi/2, \pi/2)=0.87$,
  $\chi(0.25,0.5,\pi, \pi)=1.29$,
$\chi(0.25,0.01,\pi, \pi)=1.0$.
  Concerning the consistency  of  the present schemes, both of
  the schemes are consistent with the original dif\/ferential equation (1).
  The f\/inite dif\/ference equation using the compact ADI scheme
  is consistent in sense that the local truncation error,
  $e_T = O [ \Delta x^5, \Delta t  \Delta x^2,  \Delta t \Delta x]$
  tends to zero as $\Delta t$ and $\Delta x$ tend to zero.
For Du Fort Frankel scheme equation (2.5) whose  truncation error,
$e_T  = O [ \Delta t^2 , (\Delta t/\Delta x)^2 ,  \Delta x^4]$,
the consistency condition requires the truncation error tends to zero
upon  $(\Delta t/\Delta x)^2$  approach zero as $\Delta t$ and
$\Delta x$ approach zero. For this reason,  and  a much smaller time step
than allowed by the above stability limit is implied.
This concludes that each of the f\/inite dif\/ference approximations
to the 2-D coupled Burgers' equation, the fourth-order explicit
Du Fort Frankel scheme and the compact scheme, satisf\/ies the
consistency condition.  Then, the stability of  the scheme will be the
necessary and suf\/f\/icient condition for convergence, which is true
for linear PDE's.  But, for the present nonlinear PDE's (1),
the results of the test cases will verify  the convergence, but with
higher  restricted stability limit.

\setcounter{equation}{0}

\section{Numerical experiments}

For small value of $\nu$, Burgers' equation behaves merely as hyperbolic
partial dif\/ferential equation, and the problem becomes very dif\/f\/icult
tosolve as steep shock-like wave fronts developed, as reported by Kreiss~[57].
Therefore, the present higher-order schemes are applied to solve
problems that are dominated by moderate to severe internal and boundary
gradients.

\subsection*{Problem case--1}

The f\/irst test case is the solution of 2-D unsteady coupled Burgers'
equations (1) in the domain $\{ -1 < x < 1, \ 0 < y <  \pi/6k \}$
with initial and Dirichlet boundary conditions given by the exact
steady-state solutions, that are set to form moderate to severe
internal and boundary gradients in the domain~[48], as shown in Fig.~2a:
\[
u_s(x,y)=-\frac{2}{Re}\frac{\phi_1(x,y)}{\phi(x,y)},
\qquad v_s(x,y) =-\frac{2}{Re}\frac{\phi_2(x,y)}{\phi(x,y)},
\]
where
\be \label{radwa:3.1}
\ba{l}
\ds \phi(x,y)= a_0 +a_1 x+ \left[ e^{k(x-1)}+e^{-k(x-1)}\right]\cos ky,
\vspace{3mm}\\
\ds \phi_1 (x,y)= a_1 + k\left[ e^{k(x-1)}-e^{-k(x-1)}\right]\cos ky,
\vspace{3mm}\\
\ds \phi_2 (x,y)= -k\left[ e^{k(x-1)}+e^{-k(x-1)}\right]\sin ky.
\ea
\ee
The values of the parameters $(a_0, a_1, k, Re)$ determine the type of
the gradient in the computed solutions.  We consider three cases
in this problem:

\begin{enumerate}
\item[]  {\bf case 1a.}  moderate internal gradient with
 $a_0=a_1=110.13$, $k=5$ and $Re=10$;

\item[] {\bf case 1b.}  severe internal  gradient with
$a_0=a_1=1.2962 \times 10^{13}$, $k=25$ and $Re=50$;

\item[] {\bf case 1c.}
 severe boundary gradient with
 $a_0=a_1=0.011013$, $k=5$ and $Re=10$.
\end{enumerate}

The numerical steady-state solutions of the equations (1) have been
 obtained at time $= 0.1$ for the above three cases, using the present
 schemes for dif\/ferent grid sizes.   Fig.~2b  shows the computed values
 of $u$-velocity component for two dif\/ferent grids;
 $(10\times 5)$, $(40\times 20)$, using the compact ADI scheme.
 The fourth-order accurate compact ADI scheme is
 capable of producing convergent and stable steady-state solutions with
 severe gradient even on relatively coarse grid size and large time step
 size $(\Delta t=0.01)$,
 in comparison with the fourth-order Du Fort Frankel scheme
 that required f\/iner grid $(80\times 40)$, and smaller time step size
 $(\Delta t=0.0002)$, as shown in Fig.~3.
 Moreover,  the the fourth-order Du Fort Frankel
  scheme exhibits overshoots at the steep gradients especially with
  a coarse grid.  This indicates that the fourth-order Du Fort Frankel
  scheme is unstable, and a much smaller time step and grid step
  sizes than allowed by the linear stability  condition are required.
  To test the ef\/fect of the initial conditions on the performance
  of the present schemes, initial conditions dif\/ferent from the
  ones given by equation (3.1), $u_{ij}= 1$ \& $v_{ij} = y/ y_M$,
  are used.  The results obtained by the compact ADI scheme are
  stable and convergent, while the fourth-order Du Fort Frankel scheme
  suf\/fers from overshoots at the gradients,  as  shown  in Fig.~3b.
  The computational ef\/f\/iciency of the two schemes has been tested
  by measuring the execution times, using
  PC-80486 DX2/66, necessary to obtain steady-state solutions at time $= 0.1$,
  and  by computing the two-dimensional error norms $(E_u, E_v)$
   def\/ined by:
\be \label{radwa:3.2}
E_u =\frac{1}{NM}\sum_{i=1}^{N}\sum_{j=1}^{M}|u_{ij}-u_s|,
\qquad
E_v =\frac{1}{NM}\sum_{i=1}^{N}\sum_{j=1}^{M}|v_{ij}-v_s|,
\ee
where $(u_{ij}, v_{ij})$
represent the numerical computed solutions and  $(u_s, v_s)$
the exact solutions.  Table~1 shows a comparison of the execution
times  and  the error norms of the present schemes;  the compact ADI
scheme and the fourth-order Du Fort Frankel scheme. The compact ADI scheme
computations with $\Delta t_{\max} =0.01$ required about 12.5
seconds, in comparison with 26 seconds needed by the fourth-order
Du Fort Frankel scheme with $\Delta t_{\max} =0.0002$,
for the same grid $(40\times 20)$.   Moreover, the compact ADI scheme
produces more accurate and hence more ef\/f\/icient solution, even on
a courase grid, as shown in Table~1.
This concludes  that the fourth-order accurate compact ADI scheme
is {\it twice} more economical and more accurate than the fourth-order
accurate Du Fort Frankel scheme.

\subsection*{Problem case--2}

In this test case, we consider the solution of 2-D  unsteady coupled
Burgers' equations (1), that is dominated by internal gradients, in the
domain $\{0 < x < 1, \ 0 < y < 1, \ t > 0\}$,
with the following initial and Dirichlet boundary conditions [50]:
\be \label{radwa:3.3}
u(x,y,0)=\sin({\pi x})\sin(\pi y),
\ee
\be \label{radwa:3.4}
v(x,y,0)=\{\sin(\pi x) +\sin (2\pi x)\}\{ \sin(\pi y)+\sin (2\pi y)\},
\ee
\be \label{radwa:3.5}
u(0,y,t)=u(1,y,t)=u(x,0,t)=u(x,1,t)=0,
\ee
\be \label{radwa:3.6}
v(0,y,t)=v(1,y,t)=v(x,0,t)=v(x,1,t)=0.
\ee

The computed values of the velocity components $u$ and $v$, for the case
of $Re = 1$ and at dif\/ferent times $(t = 0,\; 0.01, \; 0.05)$, on two
dif\/ferent grids $(40\times 40)$ , $(10\times 10)$ , using
the compact ADI scheme, are
shown in Fig.~4. 
Again, the fourth-order compact ADI scheme is capable of producing
stable and accurate solution with internal gradient on a coarse
grid $(10 \times 10)$, and with comparable accuracy to the solution
on a ref\/ined grid $(40\times 40)$. Tables 2-3 show a comparison of
the computed solutions using the present two schemes, with the
previous numerical results of Arminjon and Beauchamp [50], using the
f\/inite element and the method of lines.
The comparison shows that the solutions obtained by the compact ADI
scheme, and by both the f\/inite element and the method of lines
coincide to three signif\/icant digits in most cases.  The behaviour of
the solution for u and v components at larger time is also
considered.  The computed solution decreases very rapidly to zero at
time $= 0.5$, which is the same result obtained previously by
Arminjon and Beauchamp [50].  Moreover, the present compact ADI
scheme has the advantage over the other schemes that it can use a
coarse grid and time step size at least twice as larger as for the
other schemes to get convergent and accurate solution.  Fig.~5 
shows a comparison of the computed solutions for v-velocity at time
$=0.01$, using the compact ADI scheme, and using the fourth-order Du
Fort Frankel scheme that required f\/iner grid $(40\times 40)$ and
smaller time step size $(\Delta t= 10^{-6})$ to produce stable solutions.
Also, the execution times, using PC-80486 DX2/66, necessary to obtain
stable solutions at time $= 0.1$ and $Re = 1$ for the present two
schemes, are listed in Table~4.  The compact ADI scheme computations,
with $\Delta t_{\max} = 10^{-3}$, required about 85 seconds, while
the corresponding computations with the fourth-order Du Fort Frankel
scheme, with $\Delta t_{\max} =10^{-6}$, required 780 seconds, for
the same grid $(40 \times 40)$. This indicates that the present fourth-order
accurate compact ADI scheme has higher computational ef\/f\/iciency for
solving the unsteady coupled Burgers' equations with internal
gradients.

\section{Conclusion}

In conclusion, the fourth-order accurate two-point compact ADI scheme
and the fourth-order accurate Du Fort Frankel scheme are used to
solve the two-dimensional unsteady coupled Burgers' equations, for
problems that are dominated by moderate to severe internal and
boundary gradients.  The accuracy and the computational ef\/f\/iciency of
the present schemes are tested.  The compact ADI scheme is found to
be stable, ef\/f\/icient, and with better resolution of steep gradients
in comparison with the other scheme, and with the previous numerical
results of Arminjon and Beauchamp, using the f\/inite element and the
method of lines~[50].

\subsection*{Acknowledgments}
The author thank the reviewers for their valuable comments, which
improved the quality of this paper.

\newpage

\noindent
{\bf Table 1.}  Comparison of execution times and error norms for
computed solutions of 2-D coupled Burgers' equations at time $= 0.1$
\& $Re = 50$, for problem case-1.

\vspace{-2mm}

\begin{center}
\begin{tabular}{|c|c|c|c|c|c|c|}
\hline
& \multicolumn{3}{|c|}{}
& \multicolumn{3}{|c|}{}
\\[-3mm]
& \multicolumn{3}{|c|}{4-th compact scheme}
& \multicolumn{3}{|c|}{4-th Du Fort Frankel scheme}
\\
& \multicolumn{3}{|c|}{}
& \multicolumn{3}{|c|}{}
\\[-3mm]
\hline
&&&&&&\\[-3mm]
$\ba{c} \mbox{grind}\\
\mbox{points}
\ea $ & $\ba{c} \mbox{execution time}\\
\mbox{(seconds)}\ea$ & $-\log E_u$ &
$-\log E_v$ & $\ba{c} \mbox{execution time}\\
\mbox{(seconds)}\ea$  & $-\log E_u$ & $-\log E_v$\\[3mm]
\hline
&&&&&&\\[-3mm]
$10 \times 5$ & $\ba{c}\Delta t=0.01\\ \mbox{time}=1\ea$ & 12.0 & 11.4&
$\ba{c}\Delta t=0.0002\\ \mbox{time}=2.5\ea$ & 7.4 & 5.4\\[3mm]
$20 \times 10$ & $\ba{c}\Delta t=0.01\\ \mbox{time}=3\ea$ & 13.6 & 13.3&
$\ba{c}\Delta t=0.0002\\ \mbox{time}=7\ea$ & 5.4 & 4.6\\[3mm]
$40 \times 20$ & $\ba{c}\Delta t=0.01\\ \mbox{time}=12.5\ea$ & 13.8 & 13.9&
$\ba{c}\Delta t=0.0002\\ \mbox{time}=26\ea$ & 4.8 & 4.2\\[3mm]
$80 \times 40$ & $\ba{c}\Delta t=0.01\\ \mbox{time}=42\ea$ & 14.8 & 15.4&
$\ba{c}\Delta t=0.0002\\ \mbox{time}=100\ea$ & 4.5 & 4.1\\[3mm]
\hline
\end{tabular}
\end{center}

\bigskip

\noindent
{\bf Table 2.} Comparison of  the values of $u$ and $v$ computed by the Compact
ADI scheme   and Du Fort Frankel scheme  at time $= 0.01$ \& $Re =1$,
 for problem case-2.

\vspace{-2mm}

\begin{center}
\begin{tabular}{|c|ccc|ccc|}
\hline
& \multicolumn{3}{|c|}{}
& \multicolumn{3}{|c|}{}
\\[-3mm]
points & \multicolumn{3}{|c|}{compact ADI scheme}
& \multicolumn{3}{|c|}{4-th Du Fort Frankel scheme}
\\[2mm]
\hline
&\multicolumn{1}{|c|}{}&\multicolumn{1}{|c|}{}&\multicolumn{1}{|c|}{}&
\multicolumn{1}{|c|}{}&\multicolumn{1}{|c|}{}&\multicolumn{1}{|c|}{}\\[-3mm]
\multicolumn{1}{|c|}{}&
\multicolumn{1}{|c|}{$ \ba{c} N=10 \\ \Delta t =1/800\ea $}
& \multicolumn{1}{|c|}{$ \ba{c} N=20 \\ 1/800\ea $}
& \multicolumn{1}{|c|}{$ \ba{c} N=40 \\ 1/1000\ea $}
& \multicolumn{1}{|c|}{$ \ba{c} N=10 \\ \Delta t=10^{-6}\ea $}
& \multicolumn{1}{|c|}{$ \ba{c} N=20 \\ 10^{-6}\ea $}
& \multicolumn{1}{|c|}{$ \ba{c} N=40 \\ 10^{-6}\ea $}\\[3mm]
\hline
& \multicolumn{3}{|c|}{}
& \multicolumn{3}{|c|}{}
\\[-3mm]
& \multicolumn{3}{|c|}{the velocity $u$}
& \multicolumn{3}{|c|}{the velocity $u$}
\\[1mm]
(0.1,0.1) & 0.07320 & 0.07275 & 0.07273 & 0.09549 & 0.08446 & 0.07729\\
(0.2,0.8) & 0.27800 & 0.27803 & 0.27800 & 0.30231 & 0.28314 & 0.27929\\
(0.4,0.4) & 0.72292 & 0.72290 & 0.72285 & 0.72942 & 0.72215 & 0.72183\\
(0.7,0.1) & 0.20542 & 0.20506 & 0.20497 & 0.25000 & 0.21957 & 0.20978\\
(0.9,0.9) & 0.07968 & 0.07955 & 0.07953 & 0.09549 & 0.08785 & 0.08276\\[1mm]
\hline
& \multicolumn{3}{|c|}{}
& \multicolumn{3}{|c|}{}
\\[-3mm]
& \multicolumn{3}{|c|}{the velocity $v$}
& \multicolumn{3}{|c|}{the velocity $v$}
\\[1mm]
(0.1,0.1) & 0.43599 & 0.43662 & 0.43448 & 0.80425 & 0.63005 & 0.51128\\
(0.2,0.8) & \hspace{-3.3pt}-0.13444 & \hspace{-3.3pt}-0.13131
& \hspace{-3.3pt}-0.13148 & \hspace{-3.3pt}-0.29180 & \hspace{-3.3pt}-0.16784 & -0.13866\\
(0.4,0.4) & 1.65503 & 1.65869 & 1.65917 & 1.66391 & 1.65563 & 1.65177\\
(0.7,0.1) & 0.06486 & 0.06337 & 0.06417 & \hspace{-3.3pt}-0.12738
& 0.00655 & 0.04585\\
(0.9,0.9) & 0.01427 & 0.01512 & 0.01476 & 0.07771 & 0.04794 & 0.02743\\[1mm]
\hline
\end{tabular}
\end{center}

\newpage

\noindent
{\bf Table 3.} Comparison of  the values of $u$ and $v$ computed by
the method of
lines and the f\/inite element method [50]  at time $= 0.01$ \& $Re =1$,
 for problem case-2.

\vspace{-2mm}

\begin{center}
\begin{tabular}{|c|ccc|ccc|}
\hline
& \multicolumn{3}{|c|}{}
& \multicolumn{3}{|c|}{}
\\[-3mm]
points & \multicolumn{3}{|c|}{method of lines}
& \multicolumn{3}{|c|}{finite element method}
\\[2mm]
\hline
&\multicolumn{1}{|c|}{}&\multicolumn{1}{|c|}{}&\multicolumn{1}{|c|}{}&
\multicolumn{1}{|c|}{}&\multicolumn{1}{|c|}{}&\multicolumn{1}{|c|}{}\\[-3mm]
\multicolumn{1}{|c|}{}&
\multicolumn{1}{|c|}{$ \ba{c} N=10 \\ \Delta t =1/800\ea $}
& \multicolumn{1}{|c|}{$ \ba{c} N=20 \\ 1/2000\ea $}
& \multicolumn{1}{|c|}{$ \ba{c} N=40 \\ 1/4000\ea $}
& \multicolumn{1}{|c|}{$ \ba{c} N=10 \\ \Delta t=1/2000\ea $}
& \multicolumn{1}{|c|}{$ \ba{c} N=20 \\ 1/1500\ea $}
& \multicolumn{1}{|c|}{$ \ba{c} N=40 \\ 1/2000\ea $}\\[3mm]
\hline
& \multicolumn{3}{|c|}{}
& \multicolumn{3}{|c|}{}
\\[-3mm]
& \multicolumn{3}{|c|}{the velocity $u$}
& \multicolumn{3}{|c|}{the velocity $u$}
\\[1mm]
(0.1,0.1) & 0.07277 & 0.07257 & 0.07253 & 0.07279 & 0.07257 & 0.07252\\
(0.2,0.8) & 0.28887 & 0.28846 & 0.28836 & 0.28867 & 0.28842 & 0.28835\\
(0.4,0.4) & 0.72315 & 0.72205 & 0.72178 & 0.72370 & 0.72210 & 0.72179\\
(0.7,0.1) & 0.20139 & 0.20112 & 0.20106 & 0.20157 & 0.20117 & 0.20107\\
(0.9,0.9) & 0.07956 & 0.07948 & 0.07947 & 0.07951 & 0.07947 & 0.07946\\[1mm]
\hline
& \multicolumn{3}{|c|}{}
& \multicolumn{3}{|c|}{}
\\[-3mm]
& \multicolumn{3}{|c|}{the velocity $v$}
& \multicolumn{3}{|c|}{the velocity $v$}
\\[1mm]
(0.1,0.1) & 0.43857 & 0.43302 & 0.43173 & 0.44130 & 0.443357 & 0.43178\\
(0.2,0.8) & \hspace{-3.3pt}-0.13200 & \hspace{-3.3pt}-0.12387
& \hspace{-3.3pt}-0.12184 & \hspace{-3.3pt}-0.13172 & \hspace{-3.3pt}-0.12366
& \hspace{-3.3pt}-0.12180\\
(0.4,0.4) & 1.66509 & 1.65571 & 1.65335 & 1.66212 & 1.65499 & 1.65316\\
(0.7,0.1) & 0.06137 & 0.06571 & 0.06679 & 0.06306 & 0.06621 & 0.06692\\
(0.9,0.9) & 0.01459 & 0.01372 & 0.01349 & 0.01459 & 0.01367 & 0.01349\\[1mm]
\hline
\end{tabular}
\end{center}

\bigskip

\noindent
{\bf Table 4.}  Comparison of execution times for computed  solutions of
2-D  coupled Burgers' equations at time $=0.01$ \& $Re = 1.0$ for
problem case-2.

\begin{center}
\begin{tabular}{|c|c|c|}
\hline
&&\\[-3mm]
grid points & 4-th compact ADI scheme & 4-th Du Fort Frankel Scheme\\[1mm]
\hline
&&\\[-3mm]
$10 \times 10 $ & $\ba{c} \Delta t=1.25 \times 10^{-3}\\
\mbox{time}=4\;\mbox{sec} \ea $ &
$\ba{c} \Delta t=10^{-6}\\
\mbox{time}=35\;\mbox{sec} \ea $\\[3mm]
\hline
&&\\[-3mm]
$20 \times 20 $ & $\ba{c} \Delta t=1.25 \times 10^{-3}\\
\mbox{time}=14\;\mbox{sec} \ea $ &
$\ba{c} \Delta t=10^{-6}\\
\mbox{time}=240\;\mbox{sec} \ea $\\[3mm]
\hline
&&\\[-3mm]
$40 \times 40 $ & $\ba{c} \Delta t=10^{-3}\\
\mbox{time}=85\;\mbox{sec} \ea $ &
$\ba{c} \Delta t=10^{-6}\\
\mbox{time}=780\;\mbox{sec} \ea $\\[3mm]
\hline
\end{tabular}
\end{center}

\qquad

\newpage

\noindent
{\bf Figure Captions}

\strut\hfill

\noindent
{\bf Figure 1:} The couputed amplif\/ication factors (2.22) of the numerical
solution of coupled Burgers' equations (1), using the explicit fourth-order
accurate Du Fort Frankel scheme for dif\/ferent values of $|c|$ and $d$.

\newpage

\noindent
{\bf Figure 2a:} The excact steady state solutions for $u$ and $v$ velocities
of the coupled unsteady Burgers' equations for problem case-1.

\newpage

\noindent
{\bf Figure 2b:} The computed steady state solutions for $u$-velocity of the coupled
Burgers' equations for problem case-1, for two dif\/ferent grid sizes;
$(40\times 20)$, $(10\times 5)$, using the present compact ADI scheme.

\newpage

\noindent
{\bf Figure 3:}
The computed steady state solutions for $u$ and $v$ velocities of
the coupled unsteady Burgers' equations for problem case-1b, with
internal severe gradient, using compact ADI scheme, and Du Fort
Frankel scheme.

\newpage

\noindent
{\bf Figure 4:} The computed solutions for $u$ and $v$ velocities of the coupled
unsteady Burgers' equations for problem case-2 at $t=0,\; 0.01,\;0.05$, \&
$Re=1$, for two dif\/ferent grid sizes $(40\times 40)$, $(10\times 10)$,
using the present compact ADI scheme.

\newpage

\noindent
{\bf Figure 5:}
Comparison of the computed solutions of $v$-velocity component of
the coupled unsteady Burgers' equations for problem case-2, at time
$=0.01$ \& $Re=1$, using the present compact ADI scheme, and Du Fort Frankel
scheme.

\label{radwa-lp}

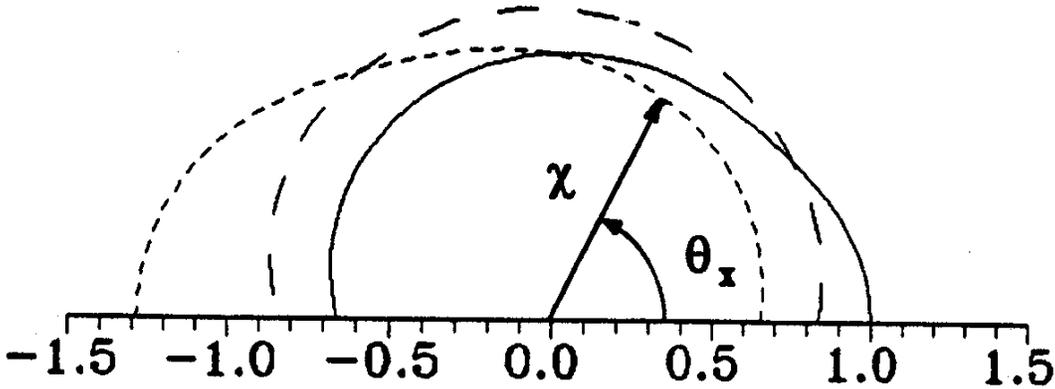
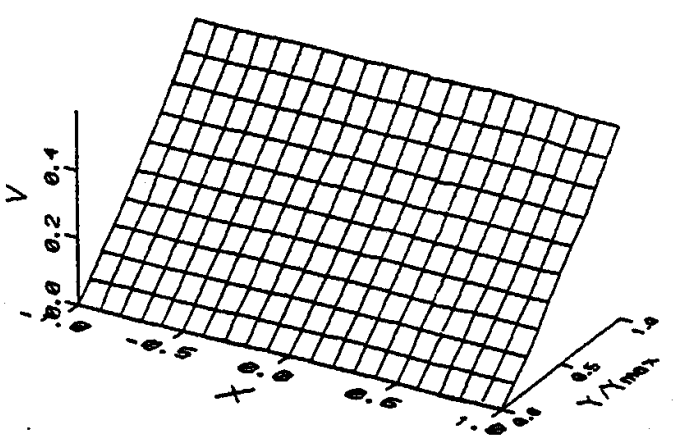
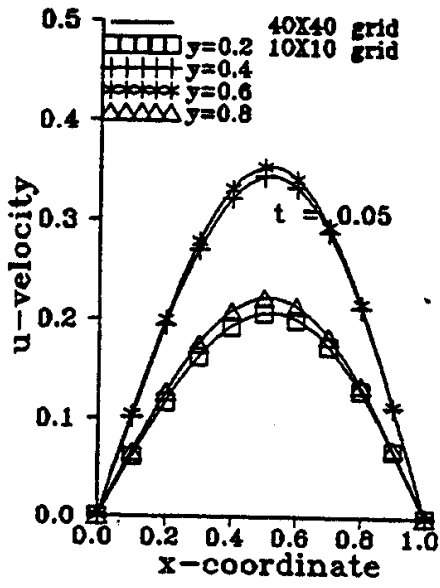
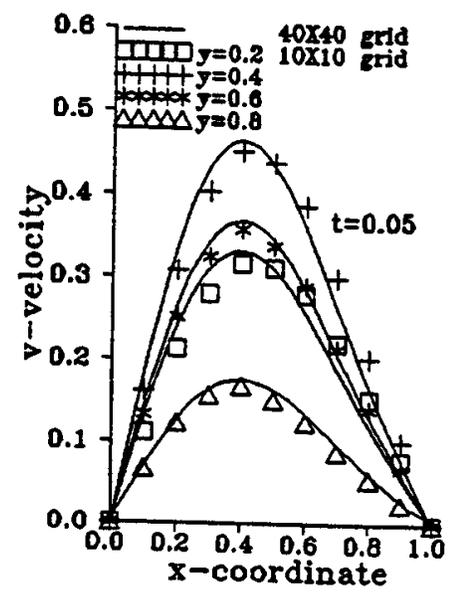
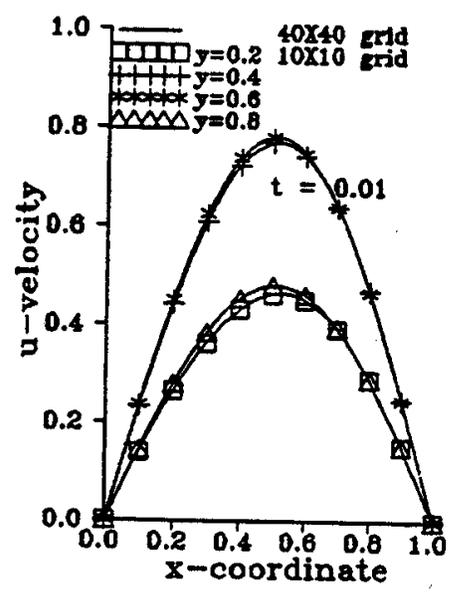
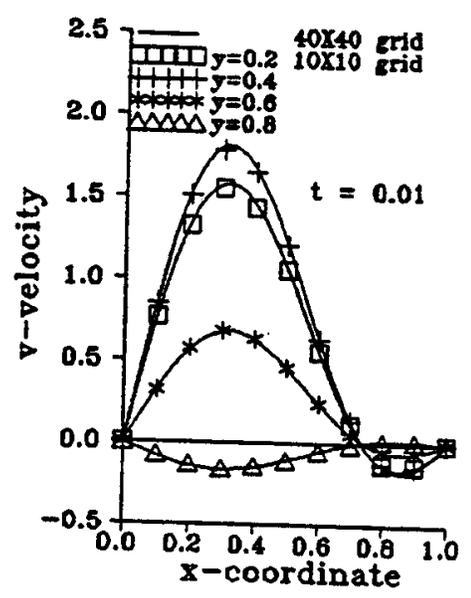
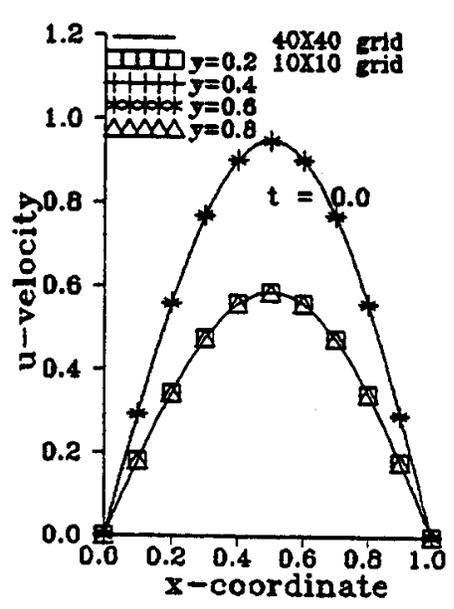
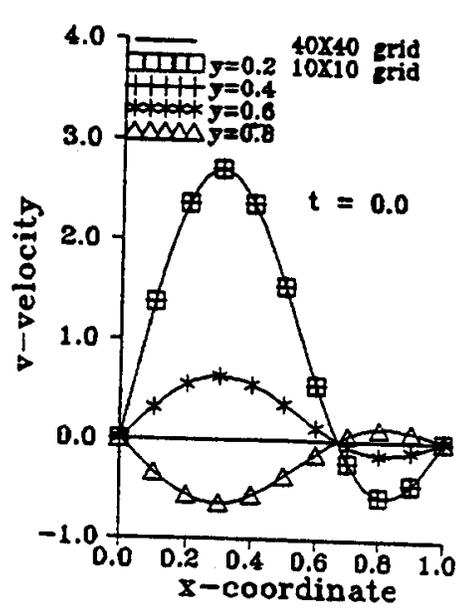
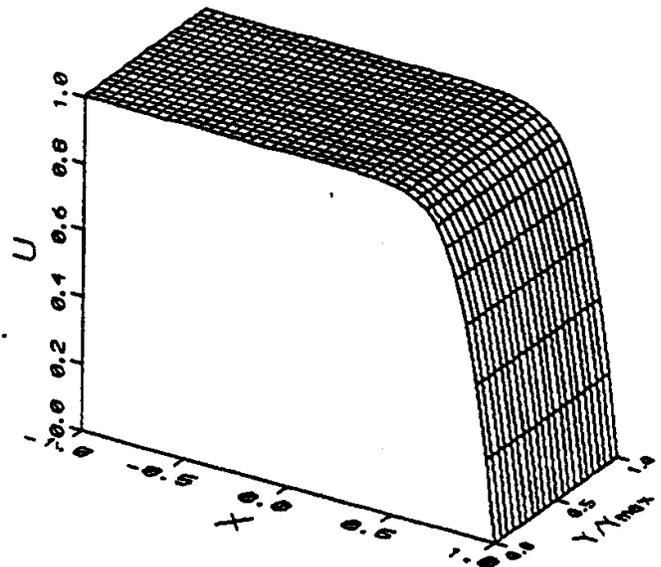
\end{document}